\documentclass[11pt]{amsart}
\usepackage{amssymb,latexsym}
\usepackage{amssymb,amsmath,ams}
\newdimen\AAdi%
\newbox\AAbo%
\def\AAk#1#2{\s_etbox\AAbo=\hbox{#2}\AAdi=\wd\AAbo\kern#1\AAdi{}}%
\def\AAr#1#2#3{\s_etbox\AAbo=\hbox{#2}\AAdi=\ht\AAbo\raise#1\AAdi\hbox{#3}}%
\font\tenmsb=msbm10 at 11pt \font\sevenmsb=msbm7 at 8pt
\font\fivemsb=msbm5 at 6pt
\newfam\msbfam
\textfont\msbfam=\tenmsb \scriptfont\msbfam=\sevenmsb
\scriptscriptfont\msbfam=\fivemsb

\newtheorem{theorem}{Theorem}[section]
\newtheorem{lemma}[theorem]{Lemma}
\theoremstyle{definition}

\newtheorem{proposition}[theorem]{Proposition}

\theoremstyle{remark}

\numberwithin{equation}{section}

\def\R{\mathbb R}
\def\S{\mathcal S}
\def\O{\mathcal O}
\def\a {\alpha}
\def\b{\beta}
\def\n{\nabla}
\begin{document}
\title  [Curvature estimates for the level sets]
{Curvature estimates for the level sets of spatial quasiconcave solutions to a class
of parabolic equations}
\author{Chuanqiang Chen}
\address{Department of Mathematics\\
         University of Science and Technology of China\\
         Hefei 230026, Anhui Province, CHINA.}
\email{cqchen@mail.ustc.edu.cn}
\author{Shujun Shi}
\address{Department of Mathematics\\
         University of Science and Technology of China\\
         Hefei 230026, Anhui Province, CHINA\\
         and
         School of Mathematical Sciences\\
         Harbin Normal University\\
         Harbin 150025, Heilongjiang Province, CHINA.}
\email{shjshi@mail.ustc.edu.cn}
\thanks{Research of the first author was supported by Grant 10871187 from the
National Natural Science Foundation of China. Research of the second
author was supported in part by the Science Research Program from
the Education Department of Heilongjiang Province (11551137).}

\begin{abstract} We prove a constant rank theorem for the second
fundamental form of the spatial convex level surfaces of solutions
to equations $u_t=F(\n^2u, \n u, u, t)$ under a structural
condition, and give a geometric lower bound of the principal
curvature of the spatial level surfaces.
\end{abstract}

\maketitle

\section{Introduction}

In this paper, we consider the convexity and principal curvature estimates of the spatial level surfaces of the spatial quasiconcave solutions to a
class of parabolic equations under some structural conditions. A
continuous function $u(x,t)$ on $\Omega \times [0,T]$ is called {\it
spatial quasiconcave} if its level sets $\{x \in \Omega| u(x,t)\ge
c\}$ are convex for each constant $c$ and any fixed $t \in [0,T]$.

The convexity of the level sets of the solutions to elliptic partial
differential equations has been studied extensively. For instance, Ahlfors \cite{AH} contains the well-known result that level curves of Green function on simply connected convex domain in the plane are the convex Jordan curves. In 1956, Shiffman \cite{Sh56} studied the minimal annulus in $\mathbb{R}^3$ whose
boundary consists of two closed convex curves in parallel planes $P_1, P_2$. He proved that the intersection of the surface with any parallel plane $P$, between $P_1$ and $P_2$, is a convex Jordan curve. In 1957, Gabriel \cite{Ga57} proved that the level sets of the Green function on a 3-dimensional bounded convex domain are strictly convex. In 1977, Lewis \cite{Le77} extended Gabriel's result to $p$-harmonic functions in higher dimensions.
Caffarelli-Spruck \cite{CS82} generalized the Lewis \cite{Le77} results to a class of semilinear elliptic partial differential equations. Motivated by the result of Caffarelli-Friedman \cite{CF85}, Korevaar \cite{Ko90} gave a new proof on the results of Gabriel and Lewis by applying the deformation process and the constant rank theorem of the second fundamental form of the convex level sets of $p$-harmonic function. A survey of this subject is given by Kawohl \cite{Ka85}. For more recent related extensions, please see the papers by Bianchini-Longinetti-Salani \cite{BLS}, Bian-Guan \cite{BG09}, Xu \cite{Xu08} and Bian-Guan-Ma-Xu \cite{BGMX}.

 There is also an extensive literature on the curvature estimates of the level sets of the solutions to elliptic partial differential equations. For 2-dimensional harmonic function and minimal surface with convex level curves, Ortel-Schneider \cite{OS}, Longinetti \cite{Lo83} and \cite{Lo87} proved that the curvature of the level curves attains its minimum on the boundary (see Talenti \cite{T83} for related
results). Longinetti also studied the precise relation between the curvature of the convex level curves and the height of 2-dimensional minimal surface in \cite{Lo87}.  Ma-Ou-Zhang \cite{MOZ09} got the Gaussian curvature estimates of the convex level sets on higher dimensional harmonic function, and Wang-Zhang \cite{WZ} got the similar curvature estimates of some quasi-linear elliptic equations under certain structure condition \cite{BLS}. Both of their test functions involved the Gaussian curvature of the boundary and the norm of the gradient on the boundary. Furthermore, for the $p$-harmonic function with strictly convex level sets, Ma-Zhang \cite{MZ} obtained that the curvature function introduced in it is concave with respect to the height of the $p$-harmornic function. For the principal curvature estimates in higher dimension, in terms of the principal curvature of the boundary and the norm of the gradient on the boundary, Chang-Ma-Yang \cite{CMY} obtained the lower bound estimates of principal curvature for the strictly convex level sets of higher dimensional harmonic functions and solutions to a class of semilinear elliptic equations under certain structure condition \cite{BLS}. Recently, in Guan-Xu \cite{GX}, they got a lower bound for the principal curvature of the level sets of solutions to a class of fully nonlinear elliptic equations in convex rings under the general structure condition \cite{BLS} via the approach of constant rank theorem.

Naturally, we hope to give a  characterization about the convexity and curvature of the level surfaces of the solutions to the corresponding parabolic equations. Borell \cite{Bo} showed the same property in \cite{Ga57} and \cite{Le77} for the solution of the corresponding heat conduction problem with zero initial data. In this paper, we will consider the following parabolic equations
\begin{equation}\label{1.1}
\frac{{\partial u}} {{\partial t}} = F(\n ^2u, \n u, u, t),
\quad~\text{in}~ \Omega  \times (0,T],
\end{equation}
where $\Omega$ is a domain in $\R^n$, and $\n ^2 u$, $\n u$ are the
spatial Hessian and spatial gradient of $u(x,t)$ respectively. Let
$\S^n$ denote the space of real symmetric $n \times n$ matrices,
$\Lambda\subset \S^n$ an open set,  $\mathbb S^{n-1}$ a unit sphere and $F=F(r,p,u,t)$  a
$C^{2,1}$ function in $\Lambda \times \R^n \times \R \times [0,T] $.
We will assume that $F$ satisfies the following conditions: there are
$\gamma_0>0$ and $c_0\in \mathbb R$,
\begin{equation}\label{1.2}
F^{\alpha\beta}:=\left(\frac{\partial F}{\partial
r_{\alpha\beta}}(r, p, u, t)\right) >0, \quad \forall \; (r, p, u,
t) \in \Lambda\times \mathbb R^{n}\times (-\gamma_0+c_0,
\gamma_0+c_0) \times [0,T],
\end{equation}
and for each $(\theta,u)\in \mathbb S^{n-1}\times \mathbb R$
fixed,
\begin{equation}\label{1.3}
F(s^2 A, s \theta, u, t) \text{ is locally concave in  } (A,s)
\text{ for each fixed }  t.
\end{equation}

Now we state our theorems.
\begin{theorem}\label{th1.1}
Suppose $u \in C^{3,1}(\Omega \times [0,T])$ is a spatial
quasiconcave solution to parabolic equation (\ref{1.1}) such that
$(\n^2u(x,t), \n u(x,t), u(x,t))\in \Lambda\times \mathbb R^n\times
(-\gamma_0+c_0, \gamma_0+c_0)$ for each $(x,t)\in \Omega \times [0,
T]$. Suppose that, $F$ satisfies conditions (\ref{1.2}) and
(\ref{1.3}),  $\n u \neq 0$ and the spatial level sets $\{x \in
\Omega |u(x,t) \ge c\}$ of $u$ are connected and locally convex for
all $c\in (-\gamma_0+c_0, \gamma_0+c_0)$ for some $\gamma_0>0$. Then
the second fundamental form of spatial level surfaces $\{x \in
\Omega |u(x,t) = c\}$ has the same constant rank for all $c\in
(-\gamma_0+c_0, \gamma_0+c_0)$. Moreover, let $l(t)$ be the minimal
rank of the second fundamental form in $\Omega$, then $l(s)
\leqslant l(t)$ for all $s \leqslant t \leqslant T$.
\end{theorem}

Inspired by \cite{GX}, we also consider to establish a geometric
lower bound for the principal curvature of the spatial level
surfaces of solutions to parabolic equation on the convex rings as
follows,
\begin{equation}\label{1.4}
\left\{\begin{array}{lcl} \frac{{\partial u}} {{\partial t}} =
F(\n ^2u, \n u, u, t)   &\text{in}& {\Omega  \times (0,T]}, \\
              u(x,0) = u_0 (x)  &\text{in}& \Omega,   \\
              u(x,t) =  0   &\text{on}&  \partial \Omega_0 \times (0,T],\\
              u(x,t) =  1   &\text{on}&  \partial \Omega_1 \times (0,T],
\end{array} \right.
\end{equation}
where $\Omega  = \Omega _0 \backslash \overline {\Omega _1 }$,
$\Omega_0$, $\Omega_1$ are two convex domains with $\overline
{\Omega _1 } \subset \Omega_0$, $F(\n ^2u_0, \n u_0, u_0, 0)>0$ and $u_0$ is quasiconcave and satisfies
\begin{equation}\label{1.5}
\left\{\begin{array}{lcl}
          u_0 =  0  & \text{on}&  \partial \Omega_0,\\
              u_0 =  1  & \text{on}&  \partial \Omega_1.
\end{array} \right.
\end{equation}
We denote $\kappa_s(x,t)$ the smallest principal curvature of
the spatial level set $\Sigma^{u(x_0,t)}=\{x \in \Omega |u(x,t) =
u(x_0,t)\}$ at $(x,t)$. For each $(x_0,t)$, set
\begin{equation}
\kappa^{u(x_0,t)}=\mathop {\inf }\limits_{x \in
\Sigma^{u(x_0,t)}}\kappa_s(x,t).
\end{equation}
We will assume that there exists $\lambda >0$, such that
\begin{equation}\label{1.7}
(F^{\alpha\beta}(\n ^2u, \n u, u, t)) \geq \lambda (\delta_{\a\b}),
\quad \forall (x,t) \in \overline{\Omega} \times [0,T].
\end{equation}

\begin{theorem}\label{th1.2}
Suppose $u \in C^{3,1}(\Omega \times [0,T])$ is a spatial
quasiconcave solution to parabolic equation (\ref{1.4}), and $F$
satisfies conditions (\ref{1.7}) and (\ref{1.3}),  $\n u \neq 0$,
then

\begin{equation}\label{1.8}
\kappa^{u(x,t)}\geq \min\{\kappa^0,\kappa^1e^{-A} \}e^{Au(x,t)}
\end{equation}
for some universal constant $A$ depending only on $ \left\| F
\right\|_{C^2 } $, $ n $, $ \lambda $, $ \mathop {\min
}\limits_{(x,t) \in \overline \Omega   \times [0,T]} \left| {\nabla
u} \right| $, $ \left\| u \right\|_{C^3 } $.
Moreover, if
$"="$ holds for some $u(x,t)\in (0, 1)$, then the $"="$ holds for all $u(x,t)\in [0, 1]$.

\end{theorem}

Theorem \ref{th1.1} and Theorem \ref{th1.2} may be looked as some parabolic versions for Theorem 1.1 in \cite{BGMX} and Theorem 1.5 in \cite{GX} respectively. The main idea to prove the main theorems in this paper can be found in the two literatures.

The rest of the paper is organized as follows. In section 2, we
prove Theorem \ref{th1.1}. In section 3, we prove Theorem
\ref{th1.2}.

\textbf{Acknowledgement} The authors would like to express sincere
gratitude to Prof. Xi-Nan Ma for his encouragement and many
suggestions in this subject.

\section{Proof of Theorem \ref{th1.1}}

Suppose $u(x,t) \in C^{3,1}(\Omega \times [0.T])$, and $u_n \ne 0$
for any fixed $(x,t)\in \Omega \times [0,T]$. It follows that the
upward inner normal direction of the spatial level sets $\{x \in
\Omega |u(x,t) = c\}$ is
\begin{eqnarray}\label{2.1}
\vec{n}= \frac{|u_n|}{|\n u|u_n}(u_1,u_2,...,u_{n-1},u_n),
\end{eqnarray}
where $\n u=(u_1,u_2,...,u_{n-1},u_n)$ is the spatial gradient of
$u$.

The second fundamental form $II$ of the spatial level surface of
function $u$ with respect to the upward normal direction (\ref{2.1})
is
\begin{equation}\label{2.2}
b_{ij} =  - \frac{|u_n|(u_n^2 u_{ij} + u_{nn}u_{i}u_j - u_nu_ju_{in}
- u_nu_iu_{jn})}{|\n u|u_n^3}.
\end{equation}
Set
\begin{equation}\label{2.3}
h_{ij} =u_n^2 u_{ij} + u_{nn}u_{i}u_j - u_nu_ju_{in} - u_nu_iu_{jn},
\end{equation}
we may write
\begin{equation}\label{2.4}
b_{ij} = - \frac{|u_n|h_{ij}}{|\n u|u_n^3}.
\end{equation}
Note that if $\Sigma ^{c,t} = \{x \in \Omega |u(x,t) = c\}$ is
locally convex, then the second fundamental form of $\Sigma^{c,t}$
is semipositive definite with respect to the upward normal direction
(\ref{2.1}). Let $a(x,t)=(a_{ij}(x,t))$ be the symmetric Weingarten
tensor of $\Sigma^{c,t}= \{x \in \Omega |u(x,t) = c \}$, then $a$ is
semipositive definite. As computed in \cite{BGMX}, if $u_n \ne 0$,
and the Weingarten tensor is
\begin{equation}\label{2.5}
a_{ij} =-\frac{|u_n|}{|\n u|{u_n}^3}\left\{h_{ij}
-\frac{u_iu_lh_{jl}}{W(1+W)u_n^2} -\frac{u_ju_lh_{il}}{W(1+W)u_n^2}
+ \frac{u_iu_ju_ku_l h_{kl}}{W^2(1+W)^2u_n^4}\right\}.
\end{equation}
With the above notations, at the point $(x,t)$ where $u_n(x,t)=|\nabla
u(x,t)|>0,\, u_i(x,t)=0$, $i=1, \cdots, n-1$, $a_{ij,k}$ is
commutative, that is, they satisfy the Codazzi property
$a_{ij,k}=a_{ik,j}, \; \forall i,j,k \le n-1$.

\subsection{Calculations on the test function}

Since Theorem \ref{th1.1} is of local feature, we may assume level
surface $\Sigma^{c,t} = \{x \in \Omega| u(x,t)=c\}$ is connected for
each $c\in (c_0-\gamma_0, c_0+\gamma_0)$. Suppose $a(x,t_0)$ attains
minimal rank $l=l(t_0)$ at some point $z_0 \in \Omega$. We may
assume $l\leqslant n-2$, otherwise there is nothing to prove.  And
we assume $u \in C^{3,1}(\Omega\times [0,T])$ and $u_n>0$ in the
rest of this paper. So there is a neighborhood $\mathcal {O}\times
(t_0-\delta, t_0+\delta]$ of $(z_0, t_0)$, such that there are $l$
"good" eigenvalues of $(a_{ij})$ which are bounded below by a
positive constant, and the other $n-1-l$ "bad" eigenvalues of
$(a_{ij})$ are very small. Denote $G$ be the index set of these
"good" eigenvalues and $B$ be the index set of "bad" eigenvalues.
And for any fixed point  $(x,t) \in \mathcal {O}\times (t_0-\delta,
t_0+\delta]$, we may express $(a_{ij})$ in a form of
(\ref{2.5}), by choosing $e_1,\cdots, e_{n-1},e_n$ such that

\begin{equation}\label{2.6}
 |\n u(x,t)|=u_n(x,t)>0\ \mbox{and}
(u_{ij}),i,j=1,..,n-1, \mbox{is diagonal at}\ (x,t).
\end{equation}

Without loss of generality we assume $ u_{11} \geq u_{22}\geq \cdots
\geq u_{n-1n-1} $. So, at $(x,t) \in \mathcal {O}\times (t_0-\delta,
t_0+\delta)$, from \eqref{2.5}, we have the matrix
$(a_{ij}),i,j=1,..,n-1,$ is also diagonal, and without loss of
generality we may assume $a_{11} \geq a_{22} \geq ... \geq a_{n-1,
n-1}$. There is a positive constant $C>0$ depending only on
$\|u\|_{C^{4}}$ and $\mathcal {O}\times (t_0-\delta, t_0+\delta]$,
such that $a_{11} \ge a_{22} \ge...\ge a_{ll} > C$ for all $(x,t)
\in \mathcal {O}\times (t_0-\delta, t_0+\delta)$. For convenience we
denote $ G = \{ 1, \cdots ,l \} $ and $ B = \{
l+1, \cdots, n-1 \} $ be the "good" and "bad" sets of indices
respectively. If there is no confusion, we also denote
\begin{eqnarray}\label{2.7}
\quad \quad \mbox{\it $G=\{a_{11},...,a_{ll}\}$ and
$B=\{a_{l+1,l+1},...,a_{n-1,n-1}\}$.}
\end{eqnarray}
Note that for any
$\delta>0$, we may choose $\mathcal {O}\times (t_0-\delta,
t_0+\delta]$ small enough such that $a_{jj} <\delta$ for all $j \in
B$ and $(x,t) \in \mathcal {O}\times (t_0-\delta, t_0+\delta]$.

For each $c$, let $a=(a_{ij})$ be the symmetric Weingarten tensor of
$\Sigma^{c,t}$. Set
\begin{eqnarray}\label{2.8}
 p(a)=\sigma_{l+1}(a_{ij}), \quad
q(a) &=& \left\{
 \begin{array}{llr}
\frac{\sigma_{l+2}(a_{ij})}{\sigma_{l+1}(a_{ij})},  &\mbox{if} \; \sigma_{l+1}(a_{ij})>0&\\
0,& \mbox{otherwise}. &
\end{array}
 \right.
\end{eqnarray}
Theorem \ref{th1.1} is equivalent to say $p(a)\equiv 0$ (defined in
(\ref{2.8}) ) in $\mathcal {O}\times (t_0-\delta, t_0]$. Since we
are dealing with general fully nonlinear equation (\ref{1.1}), as in
the case for the convexity of solutions in \cite{BG09}, there are
technical difficulties to deal with $p(a)$ alone. A key idea in
\cite{BG09} is the introduction of function $q$ as in (\ref{2.8})
and explore some crucial concavity properties of $q$. We consider
function
\begin{equation}
 \label{2.9}
 \phi(a)= p(a)+q(a),
\end{equation}
where $p$ and $q$ as in \eqref{2.8}. We will use notion $h=O(f)$ if
$|h(x,t)| \le Cf(x,t)$ for $(x,t) \in \O\times (t_0-\delta, t_0+\delta]$ with positive constant $C$ under
control.

To get around $p=0$, for $\varepsilon>0$ sufficiently small, we instead
consider
\begin{equation}\label{2.10}
\phi_\varepsilon (a)= \phi(a_\varepsilon),
\end{equation}
where $a_\varepsilon=a+\varepsilon I.$ We will also denote
$G_\varepsilon=\{a_{ii}+\varepsilon, i\in G\},$
$B_\varepsilon=\{a_{ii}+\varepsilon, i\in B\}.$

To simplify the notations, we will drop subindex $\varepsilon$
 with the understanding that all the estimates will be independent
 of $\varepsilon.$ In this setting, if we pick $\mathcal {O}\times (t_0-\delta,
t_0+\delta]$ small enough,
 there is $C>0$ independent of $\varepsilon$ such that
 \begin{equation}\label{2.11}
 \phi(a(x,t))\geq C\varepsilon, \quad \sigma_1(B)\geq C\varepsilon,
 ~\quad ~\mbox{ for all}\ (x,t) \in \mathcal {O}\times (t_0-\delta,
t_0+\delta].
 \end{equation}

In what follows, we will use $i, j, \cdots$ as indices run from $1$
to $n-1$ and use the Greek indices $\alpha, \beta, \cdots$ as
indices run from $1$ to $n$. Denote
\begin{align*}
& F^{\alpha \beta}  = \frac{{\partial F}} {{\partial u_{\alpha
\beta} }} , F^{p_\alpha } = \frac{{\partial F}} {{\partial u_\alpha
}} , F^u = \frac{{\partial F}}{{\partial u}} ,F^t  = \frac{{\partial F}}{{\partial t}} ,\\
&F^{\alpha \beta,\gamma \eta}  = \frac{{\partial ^2 F}} {{\partial
u_{\alpha \beta}\partial u_{\gamma \eta} }},  F^{\alpha
\beta,p_\gamma } = \frac{{\partial ^2 F}} {{\partial u_{\alpha
\beta}\partial u_\gamma }} , F^{\alpha \beta,u} = \frac{{\partial
^2 F}}{{\partial u_{\alpha \beta}\partial u}},\\
& F^{p_\alpha p_\beta }  = \frac{{\partial ^2 F}} {{\partial
u_\alpha\partial u_\beta}} ,  F^{p_\alpha, u}  = \frac{{\partial ^2
F}} {{\partial u_\alpha\partial u}} ,  F^{u,u}  = \frac{{\partial ^2
F}} {{\partial u}^2}.
\end{align*}
We also denote
\begin{equation}
\mathcal H_{\phi} = \sum_{i,j\in B}|\nabla a_{ij}|+\phi.
\end{equation}

\begin{lemma}\label{lem2.1}
For any fixed $(x,t) \in \mathcal {O}\times (t_0-\delta,
t_0+\delta]$, with the coordinate chosen as in (\ref{2.6}) and
(\ref{2.7}),
\begin{equation}\label{2.13}
\phi_t=-u_n^{-3}\sum_{j \in B} \left[ \sigma_l(G) +
\frac{{\sigma}^2_1(B|j)-{\sigma}_2(B|j)}
{{\sigma}^2_1(B)}\right][u_n^2u_{jjt}-2u_nu_{jn}u_{jt} ]+O(\mathcal
H_{\phi})
\end{equation}
and
\begin{eqnarray}\label{3.9}
&&\sum_{\alpha,\beta=1}^{n}F^{\alpha\beta}\phi_{\alpha\beta} \notag
\\
&=&u_n^{-3}\sum_{j\in B}\left[ \sigma_l(G) +
\frac{{\sigma}^2_1(B|j)-{\sigma}_2(B|j)} {{\sigma}^2_1(B)}\right]
[-u_n^2\sum_{\alpha,\beta=1}^n F^{\alpha\beta}u_{\alpha\beta j
j}+2u_nu_{nj}\sum_{\alpha,\beta=1}^nF^{\alpha\beta}u_{\alpha\beta j}
\notag \\
&& \qquad\qquad\qquad\qquad\qquad\qquad\qquad\qquad\quad
+4u_nu_{nj}\sum_{\alpha,\beta=1}^nF^{\alpha\beta}u_{\alpha\beta
j}-6u_{nj}^2\sum_{\alpha,\beta=1}^{n}F^{\alpha\beta}u_{\alpha\beta
}]  \notag  \\
&&+2u_n^{-3}\sum_{j\in B,i\in G}\left[ \sigma_l(G) +
\frac{{\sigma}^2_1(B|j)-{\sigma}_2(B|j)}
{{\sigma}^2_1(B)}\right]\sum_{\alpha,\beta=1}^{n}
F^{\alpha\beta}\frac{1}{u_{ii}}[u_nu_{ij\alpha}-2u_{i\a}u_{jn}][u_nu_{ij\b}-2u_{i\b}u_{jn}] \notag \\
&&-\frac{1}{{\sigma}^3_1(B)}\sum_{\alpha,\beta=1}^n\sum_{i\in
B}F^{\alpha\beta}[{\sigma}_1(B)a_{ii,\alpha}-a_{ii}\sum_{j\in
B}a_{jj,\alpha}][{\sigma}_1(B)a_{ii,\beta}-a_{ii}\sum_{j\in
B}a_{jj,\beta}]\nonumber\\
&& -\frac{1}{{\sigma}_1(B)}\sum_{\alpha,\beta=1}^n\sum_{i\neq j\in
B}F^{\alpha\beta}a_{ij,\alpha}a_{ij,\beta}+O(\mathcal H_{\phi}).
\notag
\end{eqnarray}
\end{lemma}

{\bf Proof:} For any fixed point $(x,t) \in \mathcal
{O}\times (t_0-\delta, t_0+\delta]$, choose a coordinate system as
in (\ref{2.6}) so that $|\n u|=  u_n >0$ and the matrix
$(a_{ij}(x,t))$ is diagonal for $1 \le i,j\le n-1$ and nonnegative.
From the definition of $\phi$,
\begin{eqnarray}\label{2.14}
a_{jj}=-\frac{h_{jj}}{u^3_n}=-\frac{u_{jj}}{u_n}=O(\mathcal
H_{\phi}) , \forall j\in B,
\end{eqnarray}
and
\begin{eqnarray}\label{2.15}
\phi_t &=& \sum_{j \in B} \left[ \sigma_l(G) +
\frac{{\sigma}^2_1(B|j)-{\sigma}_2(B|j)} {{\sigma}^2_1(B)}\right]
a_{jj,t}+O(\mathcal H_{\phi})  \notag  \\
&=& -u_n^{-3}\sum_{j \in B} \left[ \sigma_l(G) +
\frac{{\sigma}^2_1(B|j)-{\sigma}_2(B|j)}
{{\sigma}^2_1(B)}\right][u_n^2u_{jj,t}-2u_nu_{jn}u_{jt} ]+O(\mathcal
H_{\phi})
\end{eqnarray}
Using relationship (\ref{2.14}), we have
\begin{eqnarray}\label{2.16}
\phi_{\alpha \beta}&=&\sum_{j\in B}\left[ \sigma_l(G) +
\frac{{\sigma}^2_1(B|j)-{\sigma}_2(B|j)}
{{\sigma}^2_1(B)}\right]\Big[a_{jj,\alpha \beta}-2\sum_{i\in
G}\frac{a_{ij,\alpha}a_{ij,\beta}}{a_{ii}}\Big]
\nonumber\\
&&-\frac{1}{{\sigma}^3_1(B)}\sum_{i\in
B}\Big[{\sigma}_1(B)a_{ii,\alpha}-a_{ii}\sum_{j\in
B}a_{jj,\alpha}\Big]\Big[{\sigma}_1(B)a_{ii,\beta}-a_{ii}\sum_{j\in
B}a_{jj,\beta}\Big]\nonumber\\
&&-\frac{1}{{\sigma}_1(B)}\sum_{i\neq j\in
B}a_{ij,\alpha}a_{ij,\beta}+O(\mathcal H_{\phi}).
\end{eqnarray}

So far, we have followed standard calculations as in \cite{ GM03,
BGMX, BG09}. Since $u_k=0$ for $k=1,\cdots, n-1$, from
(\ref{2.5}),
\begin{eqnarray}\label{2.17}
u_n u_{ij\a}=-u_n^2 a_{ij,\a}+u_{nj}u_{i\a}
+u_{ni}u_{j\a}+u_{n\a}u_{ij}, \quad \forall\; i,j \le
n-1,\end{eqnarray}
and for each $j\in B$,
\begin{eqnarray}\label{2.18}
a_{jj,\alpha \beta}&=&-\frac{1}{u_n^3}h_{jj,\alpha \beta}+O(\mathcal
H_{\phi})  \notag  \\
&=&-\frac{1}{u_n^3}[u_n^2u_{jj\alpha\beta}+2u_{nn}u_{j\a}u_{j\b}+2u_{n\a}u_{nj}u_{j\b}+2u_{n\b}u_{nj}u_{j\a}
\notag   \\
&& \qquad \quad  -2u_nu_{nj}u_{\a\b
j}-2u_nu_{j\a}u_{nj\b}-2u_nu_{j\b}u_{nj\a}]+O(\mathcal H_{\phi}).
\end{eqnarray}
Hence, for $j\in B$,
\begin{eqnarray}\label{2.19}
\sum_{\alpha,\beta=1}^n F^{\alpha\beta}a_{jj,\alpha \beta} =&&
\sum_{\alpha,\beta=1}^n
\frac{F^{\alpha\beta}}{u_n^3}[-u_n^2u_{\alpha\beta
jj}-4u_{n\a}u_{nj}u_{j\b} +4u_nu_{j\a}u_{nj\b} \nonumber \\
&& \qquad \qquad \quad + 2u_nu_{nj}u_{\a\b j}-2u_{nn}u_{j\a}u_{j\b}]
+O(\mathcal H_{\phi}).
\end{eqnarray}
Using the fact that
$\sum_{\alpha=1}^n
F^{\alpha n}u_{n\a}=(\sum_{\alpha,\beta=1}^n
-\sum_{\b
=1}^{n-1}\sum_{\a=1}^n)F^{\alpha\beta}u_{\a\b}$, $\forall j\in B$,
\begin{eqnarray*}\label{3.24}
\sum_{\alpha,\beta=1}^n
F^{\alpha\beta}u_{n\a}u_{j\b}=u_{nj}(\sum_{\alpha,\beta=1}^n
- \sum_{\b
=1}^{n-1}\sum_{\a=1}^n)F^{\alpha\beta}u_{\a\b} +O(\mathcal H_{\phi}),
\end{eqnarray*}
\begin{eqnarray*}\label{3.25}
\sum_{\alpha,\beta=1}^n
F^{\alpha\beta}u_{j\a}u_{nj\b}
=u_{nj}(\sum_{\alpha,\beta=1}^n
-\sum_{\a=1}^{n-1}\sum_{\b=1}^n)F^{\alpha\beta}u_{\a\b
j}+O(\mathcal H_{\phi}),
\end{eqnarray*}
and
\begin{eqnarray*}\label{3.26}
&&-2u_{nn}\sum_{\alpha,\beta=1}^n
F^{\alpha\beta}u_{j\a}u_{j\b}=-2u_{nn}F^{nn}u_{nj}^2+O(\mathcal H_{\phi})\\
&=&-2u_{nj}^2\sum_{\alpha,\beta=1}^n
F^{\alpha\beta}u_{\a\b}+4u_{nj}^2\sum_{\a=1}^{n-1}F^{\a
n}u_{n\a}+2u_{nj}^2 \sum_{\alpha,\beta=1}^{n-1}
F^{\alpha\beta}u_{\a\b}+O(\mathcal H_{\phi})\nonumber.
\end{eqnarray*}
Put above to (\ref{2.19}),
\begin{eqnarray}\label{2.20}
&&\sum_{j\in B}\sum_{\alpha,\beta=1}^n F^{\alpha\beta}
u_n^3a_{jj,\alpha\beta}\nonumber\\
&=&-u_n^2\sum_{j\in B}\sum_{\alpha,\beta=1}^n
F^{\alpha\beta}u_{\alpha\beta jj}+6u_n\sum_{j\in
B}u_{nj}\sum_{\alpha,\beta=1}^n F^{\alpha\beta}u_{\a\b j}\nonumber\\
&&-6\sum_{j\in B}u_{nj}^2\sum_{\alpha,\beta=1}^n
F^{\alpha\beta}u_{\a\b}-4u_n\sum_{j\in
B}u_{nj}\sum_{\a=1}^{n-1}\sum_{\b=1}^nF^{\alpha\beta}u_{\a\b j}\nonumber\\
&&+8\sum_{j\in B}u_{nj}^2\sum_{\a=1}^{n-1}F^{\a
n}u_{n\a}+6\sum_{j\in B}u_{nj}^2 \sum_{\alpha,\beta=1}^{n-1}
F^{\alpha\beta}u_{\a\b}+O(\mathcal H_{\phi}).
\end{eqnarray}
By (\ref{2.17}), for $j\in B$,
\begin{eqnarray}\label{2.21}
&& u_n\sum_{\a=1}^{n-1}\sum_{\b=1}^nF^{\alpha\beta}u_{\a\b
j}=u_n\sum_{\a=1}^{n}\bigg(\sum_{i \in
B}F^{\alpha i}u_{ij\a}+\sum_{i \in G}F^{\alpha
i}u_{ij\a}\bigg)\nonumber\\
&=& \sum_{\a=1}^{n} \sum_{i \in G}F^{\alpha i}
(-u_n^2 a_{ij,\a}+u_{i\a}u_{jn} + u_{j\a}u_{in})\nonumber\\
&&+
\sum_{\a=1}^{n}\sum_{i \in B}F^{\alpha
i}(u_{i\a}u_{jn} + u_{j\a}u_{in})+O(\mathcal H_{\phi})\nonumber\\
&=&-u_n^2\sum_{\a=1}^{n} \sum_{i \in G}F^{\alpha i}
a_{ij,\a}+u_{nj}\sum_{i\in G}F^{ii}u_{ii}
+2u_{nj}(\sum_{i=1}^{n-1}F^{ni}u_{ni})+O(\mathcal H_{\phi}).
\end{eqnarray}
(\ref{2.20}) and (\ref{2.21}) yield
\begin{eqnarray}\label{3.29}
\sum_{\alpha,\beta=1}^n F^{\alpha\beta} u_n^3a_{jj,\alpha\beta}
&=&-u_n^2\sum_{\alpha,\beta=1}^n F^{\alpha\beta}u_{\alpha\beta
jj}+2u_n u_{nj}\sum_{\alpha,\beta=1}^n F^{\alpha\beta}u_{\a\b j} \notag  \\
&&+4u_n u_{nj}\sum_{\alpha,\beta=1}^n F^{\alpha\beta}u_{\a\b j}-6
u_{nj}^2\sum_{\alpha,\beta=1}^n F^{\alpha\beta}u_{\a\b} \notag \\
&&+4u_n^2 u_{nj}\sum_{\a=1}^{n} \sum_{i \in G}F^{\alpha i}
a_{ij,\a}+2 u_{nj}^2\sum_{i\in G}F^{ii}u_{ii} +O(\mathcal H_{\phi}).
\end{eqnarray}
So,
\begin{eqnarray}\label{3.9}
&&\sum_{\alpha,\beta=1}^{n}F^{\alpha\beta}\phi_{\alpha\beta} \notag
\\
&=&u_n^{-3}\sum_{j\in B}\left[ \sigma_l(G) +
\frac{{\sigma}^2_1(B|j)-{\sigma}_2(B|j)} {{\sigma}^2_1(B)}\right]
[-u_n^2\sum_{\alpha,\beta=1}^n F^{\alpha\beta}u_{\alpha\beta j
j}+2u_nu_{nj}\sum_{\alpha,\beta=1}^nF^{\alpha\beta}u_{\alpha\beta j}
\notag \\
&& \qquad\qquad\qquad\qquad\qquad\qquad\qquad\qquad\quad
+4u_nu_{nj}\sum_{\alpha,\beta=1}^nF^{\alpha\beta}u_{\alpha\beta
j}-6u_{nj}^2\sum_{\alpha,\beta=1}^{n}F^{\alpha\beta}u_{\alpha\beta
}]  \notag  \\
&&-2\sum_{j\in B,i\in G}\left[ \sigma_l(G) +
\frac{{\sigma}^2_1(B|j)-{\sigma}_2(B|j)}
{{\sigma}^2_1(B)}\right]\left[\sum_{\alpha,\beta=1}^n
F^{\alpha\beta}\frac{a_{ij,\alpha}a_{ij,\beta}}{a_{ii}}
-2\frac{u_{nj}}{u_n}\sum_{\a=1}^{n} F^{\alpha i}
a_{ij,\a}-\frac{u_{nj}^2}{u_n^3}F^{ii}u_{ii}\right]  \notag \\
&&-\frac{1}{{\sigma}^3_1(B)}\sum_{\alpha,\beta=1}^n\sum_{i\in
B}F^{\alpha\beta}[{\sigma}_1(B)a_{ii,\alpha}-a_{ii}\sum_{j\in
B}a_{jj,\alpha}][{\sigma}_1(B)a_{ii,\beta}-a_{ii}\sum_{j\in
B}a_{jj,\beta}]\nonumber\\
&& -\frac{1}{{\sigma}_1(B)}\sum_{\alpha,\beta=1}^n\sum_{i\neq j\in
B}F^{\alpha\beta}a_{ij,\alpha}a_{ij,\beta}+O(\mathcal H_{\phi}).
\notag
\end{eqnarray}
In fact, for any $i \in G, j \in B$,

\begin{eqnarray}\label{2.23}
&&\sum_{\alpha,\beta=1}^n
F^{\alpha\beta}\frac{a_{ij,\alpha}a_{ij,\beta}}{a_{ii}}
-2\frac{u_{nj}}{u_n}\sum_{\a=1}^{n} F^{\alpha i}
a_{ij,\a}-\frac{u_{nj}^2}{u_n^3}F^{ii}u_{ii} \notag \\
&=&-\frac{1}{u_n^3}[\sum_{\alpha,\beta=1}^n
F^{\alpha\beta}\frac{h_{ij,\alpha}h_{ij,\beta}}{h_{ii}}
-2\frac{u_{nj}}{u_n}\sum_{\a=1}^{n} F^{\alpha i}
h_{ij,\a}+u_{nj}^2F^{ii}u_{ii}]  \notag \\
&=&-\frac{1}{u_n^3}\left\{\sum_{\alpha,\beta=1}^{n-1}
F^{\alpha\beta}\frac{1}{u_n^2u_{ii}}[u_n^2u_{ij\alpha}-u_nu_{i\a}u_{jn}][u_n^2u_{ij\b}-u_nu_{i\b}u_{jn}]
\right. \notag  \\
&& \qquad \qquad  +2\sum_{\alpha=1}^{n-1}
F^{\alpha n}\frac{1}{u_n^2u_{ii}}[u_n^2u_{ij,\alpha}-u_nu_{i\a}u_{jn}][u_n^2u_{ijn}-2u_nu_{in}u_{jn}]  \notag  \\
&& \qquad \qquad + F^{nn}\frac{1}{u_n^2u_{ii}}[u_n^2u_{ijn}-2u_nu_{in}u_{jn}][u_n^2u_{ijn}-2u_nu_{in}u_{jn}]  \notag  \\
&& \qquad \qquad -2\sum_{\a=1}^{n-1} F^{\alpha
i}\frac{1}{u_n^2u_{ii}}[u_n^2u_{ij\a}-2u_nu_{i\a}u_{jn}][u_nu_{ii}u_{nj}] \notag \\
&& \qquad \qquad -2 F^{n i}\frac{1}{u_n^2u_{ii}} [u_n^2u_{ijn}-2u_nu_{in}u_{jn}][u_nu_{ii}u_{nj}]  \notag \\
&& \qquad \qquad
\left.+F^{ii}\frac{1}{u_n^2u_{ii}}(u_nu_{ii}u_{nj})^2\right\}
\notag \\
&=&-\frac{1}{u_n^3}\sum_{\alpha,\beta=1}^{n}
F^{\alpha\beta}\frac{1}{u_n^2u_{ii}}[u_n^2u_{ij\alpha}-2u_nu_{i\a}u_{jn}][u_n^2u_{ij\b}-2u_nu_{i\b}u_{jn}].
\end{eqnarray}
Obviously, we can get
\begin{eqnarray} \label{2.24}
\sum_{\alpha,\beta=1}^n
F^{\alpha\beta}\frac{a_{ij,\alpha}a_{ij,\beta}}{a_{ii}}
-2\frac{u_{nj}}{u_n}\sum_{\a=1}^{n} F^{\alpha i}
a_{ij,\a}-\frac{u_{nj}^2}{u_n^3}F^{ii}u_{ii}\geq 0,
\end{eqnarray}
this is the Claim in \cite{BGMX}.

From the above formulas, Lemma \ref{2.1} holds. \qed

\subsection{Proof of Theorem \ref{th1.1}}

We start this section with a discussion on structure condition
\eqref{1.3}. For any function $F(r,p,u,t)$,  denote
$F^{\alpha\beta}=\frac{\partial F}{\partial r_{\alpha\beta}},
F^{p_l}=\frac{\partial F}{\partial u_l}, \cdots$ as partial
derivatives of $F$ with respect to corresponding arguments.
\begin{lemma}\label{lem2.2}
If $F$ satisfies condition \eqref{1.3}, then
\begin{eqnarray}\label{2.25}
Q(V,V)&=&F^{\alpha\beta,
\gamma \eta}X_{\alpha\beta}X_{\gamma \eta}+2F^{\alpha\beta,p_l}\theta_lX_{\alpha\beta}Y+F^{p_k, p_l}\theta_k\theta_l Y^2\nonumber\\
&&+4s ^{-1}F^{\alpha\beta}X_{\alpha\beta}Y-6F^{\alpha\beta}A_{\alpha\beta}Y^2\nonumber\\
&\leqslant& 0,
\end{eqnarray}
for every $(X_{\alpha\beta}, Y) = ((s^2\widetilde{X}_{\alpha\beta} +
2s A_{\alpha\beta}\widetilde{Y}),\widetilde{Y})$, with any
$\widetilde{V}=((\widetilde{X}_{\alpha\beta}),\widetilde{Y}) \in
\S^{n}\times \mathbb R$, where $F^{\alpha\beta,
rs},F^{\alpha\beta,u_l}, etc.$ are evaluated at $(s^{2}A, s\theta,
u, t)$, and the Einstein summation convention is used.
\end{lemma}
{\bf Proof:}
 Denoting $\widetilde{F}(A,s)=F(s^2 A, s\theta, u, t),$ condition
(\ref{1.3}) implies that $\widetilde{F}(A, s)$ is locally concave,
that is,
\begin{eqnarray}\label{2.26}
\widetilde{F}^{\alpha\beta, \gamma \eta} \widetilde{X}
_{\alpha\beta} \widetilde{X}_{\gamma \eta} + 2
\widetilde{F}^{\alpha\beta,s} \widetilde{X}_{\alpha\beta}
\widetilde{Y} + \widetilde{F}^{s,s} \widetilde{Y}^2 \leq 0,
\end{eqnarray}
for any $\widetilde{V}=((\widetilde{X}_{\alpha\beta}),\widetilde{Y})
\in \S^{n}\times \mathbb R$.

At $(A,s)$,
\begin{eqnarray*}
&& \widetilde{F}^{\alpha \beta ,\gamma \eta }  = F^{\alpha \beta
,\gamma \eta } s^2  \cdot s^2,  \\
&&\widetilde{F}^{\alpha \beta ,s}  = F^{\alpha \beta ,\gamma \eta }
s^2 \cdot 2sA_{\gamma \eta }  + F^{\alpha \beta ,p_l } s^2  \cdot
\theta _l  + F^{\alpha \beta } 2s,  \\
&&\widetilde{F}^{s,s}  = F^{\alpha \beta ,\gamma \eta } 2sA_{\alpha
\beta }  \cdot 2sA_{\gamma \eta }  + 2F^{\alpha \beta ,p_l }
2sA_{\alpha \beta }  \cdot \theta _l  + F^{p_k ,p_l } \theta _k
\cdot \theta _l  + F^{\alpha \beta } 2A_{\alpha \beta }.
\end{eqnarray*}
Set
\begin{eqnarray}
&& X_{\alpha \beta }  = s^2 \widetilde{X}_{\alpha \beta }  +
2sA_{\alpha \beta } \widetilde{Y},  \\
&&Y=\widetilde{Y},
\end{eqnarray}
so (\ref{2.26}) is equivalent to
\begin{eqnarray}\label{3.61}
&&F^{\alpha\beta,\gamma \eta}X_{\alpha\beta}X_{\gamma \eta}+2F^{\alpha\beta,p_l}\theta_lX_{\alpha\beta}Y+F^{u_k, p_l}\theta_k\theta_l Y^2\nonumber\\
&&+4s ^{-1}F^{\alpha\beta}s^2 \widetilde{X}_{\alpha\beta}\widetilde{Y}+2F^{\alpha\beta}A_{\alpha\beta}\widetilde{Y}^2\nonumber\\
&=&F^{\alpha\beta,\gamma \eta}X_{\alpha\beta}X_{\gamma \eta}+2F^{\alpha\beta,p_l}\theta_lX_{\alpha\beta}Y+F^{p_k, p_l}\theta_k\theta_l Y^2\nonumber\\
&&+4s ^{-1}F^{\alpha\beta}X_{\alpha\beta}Y-6F^{\alpha\beta}A_{\alpha\beta}Y^2\nonumber\\
&\leqslant& 0.\notag
\end{eqnarray}
Therefore, (\ref{2.25}) follows from above, and Lemma \ref{2.2} holds. \qed

Theorem~\ref{th1.1} is a direct consequence of the following
proposition and the strong maximum principle.

\begin{proposition}\label{thm2.3}
Suppose that the function $F, u$ satisfy assumptions in Theorem
\ref{th1.1}. If the second fundamental form $b_{ij}$ of
$\Sigma^{c,t_0}$ attains minimum rank $l=l(t_0)$ at certain point
$x_0\in \Omega$, then there exist a neighborhood $\mathcal {O}
\times(t_0-\delta_0,t_0+\delta_0]$ of $(x_0,t_0)$ and a positive
constant $C$ independent of $\phi$ (defined in (\ref{2.9})), such
that
\begin{equation}\label{2.29}
\sum_{\alpha,\beta=1}^{n}F^{\alpha\beta}\phi_{\alpha\beta}(x,t)-\phi_t\leq
C(\phi+|\nabla \phi|),~~\forall ~(x,t)\in \mathcal {O}
\times(t_0-\delta_0,t_0+\delta_0].
\end{equation}
\end{proposition}

{\bf Proof:}  Let $u\in C^{3,1}(\Omega
\times [0,T])$ be a spatial quasiconcave solution of equation
(\ref{1.1}) and $(u_{ij}) \in \S^n.$ Let $l=l(t_0)$ be the minimum
rank of the second fundamental forms $h_{ij}$ of $\Sigma^{c,t_0}$
($l \in \{0,1,...,n-1\}$) for every $c$ in $(-\gamma_0+c_0,
\gamma_0+c_0)$, suppose the minimum rank $l$ arrives at point $x_0
\in \Sigma^{c,t_0}$. We work on a small open neighborhood $\mathcal
{O} \times(t_0-\delta_0,t_0+\delta_0]$ of $(x_0,t_0)$. We may assume
$l \le n-2$. Lemma \ref{2.1} implies $\phi \in C^{1,1}(\mathcal {O}
\times(t_0-\delta_0,t_0+\delta_0]),$ $ \phi(x,t)\geq 0, ~\qquad~
\phi(x_0,t_0)=0 $. For $\epsilon>0$ sufficient small, let
$\phi_{\epsilon}$ defined as in (\ref{2.9}) and (\ref{2.10}), we
need to verify (\ref{2.29}) for each point $(x,t) \in \mathcal {O}
\times(t_0-\delta_0,t_0+\delta_0]$. For each fixed $(x,t)$, choose a
local coordinate $e_1,\cdots, e_{n-1}, e_n$ such that (\ref{2.6}) and
(\ref{2.7}) are satisfied. We want to establish differential
inequality (\ref{2.29}) for $\phi_\varepsilon$ defined in
(\ref{2.10}) with constant $C$ independent of $\varepsilon$. Note that we will omit the subindex $\varepsilon$ with the
understanding that all the estimates are independent of
$\varepsilon$.

By Lemma~\ref{lem2.1},
\begin{eqnarray}\label{2.30}
&&\sum_{\alpha,\beta=1}^nF^{\alpha\beta}\phi_{\alpha\beta}-\phi_t\nonumber\\
&\le& -u_n^{-3}\sum_{j\in
B}\left[\sigma_l(G)+\frac{{\sigma}^2_1(B|j)-{\sigma}_2(B|j)}{{\sigma}^2_1(B)}\right]\Big[
u_n^2(\sum_{\alpha,\beta=1}^nF^{\alpha\beta}u_{jj\alpha\beta}-u_{jjt})
\nonumber\\
&&-2u_nu_{jn}(\sum_{\alpha,\beta=1}^n
F^{\alpha\beta}u_{j\alpha\beta}-u_{jt})-4u_nu_{jn}\sum_{\alpha,\beta=1}^n
F^{\alpha\beta}u_{j\alpha\beta} +6u_{jn}^2\sum_{\alpha,\beta=1}^n
F^{\alpha\beta}u_{\alpha\beta}\Big]\nonumber\\
&&-\frac{1}{{\sigma}^3_1(B)}\sum_{\alpha,\beta=1}^n\sum_{i\in
B}F^{\alpha\beta}[{\sigma}_1(B)a_{ii,\alpha}-a_{ii}\sum_{j\in
B}a_{jj,\alpha}][{\sigma}_1(B)a_{ii,\beta}-a_{ii}\sum_{j\in
B}a_{jj,\beta}]\nonumber\\
&&-\frac{1}{{\sigma}_1(B)}\sum_{\alpha,\beta=1}^n\sum_{i\neq j,
i,j\in B}F^{\alpha\beta}a_{ij,\alpha}a_{ij,\beta}+O(\mathcal H_{\phi}).
\end{eqnarray}

For each $j\in B$, differentiating equation (\ref{1.1}) in $e_j$ direction at
$x$,
\begin{equation}\label{2.31}
u_{jt}=\sum_{\alpha,\beta=1}^nF^{\alpha\beta}u_{\alpha\beta
j}+F^{u_n}u_{jn}+O(\mathcal H_{\phi}),
\end{equation}
and
\begin{eqnarray}\label{2.32}
u_{jjt}&=&\sum_{\alpha,\beta=1}^nF^{\alpha\beta}u_{\alpha\beta
jj}+\sum_{\alpha,\beta,r,s=1}^nF^{\alpha\beta, rs}u_{\alpha\beta
j}u_{rsj} +2\sum_{\alpha,\beta,l=1}^nF^{\alpha\beta,
u_l}u_{\alpha\beta j}u_{lj}\nonumber\\
&&+2\sum_{\alpha,\beta=1}^nF^{\alpha\beta, u}u_{j\alpha\beta}u_j
+\sum_{l,s=1}^nF^{u_l, u_s}u_{lj}u_{sj}-2\sum_{l=1}^nF^{u_l,
u}u_{lj}u_j\nonumber\\
&&+F^{u,u}u_j^2+\sum_{l=1}^nF^{u_l}u_{ljj}+F^{u}u_{jj}.
\end{eqnarray}
It follows from (\ref{2.17}) that, at $(x,t)$
\begin{eqnarray}\label{2.33}
\sum_{\alpha,\beta=1}^nF^{\alpha\beta}u_{\alpha\beta
jj}-u_{jjt}&=&-\sum_{\alpha,\beta,r,s=1}^nF^{\alpha\beta,
rs}u_{\alpha\beta j}u_{rsj}
-2\sum_{\alpha,\beta=1}^nF^{\alpha\beta, u_n}u_{j\alpha\beta}u_{nj}\nonumber\\
&&-F^{u_n, u_n}u^2_{jn}-2\frac{F^{u_n}}{u_n}u^2_{jn}+O(\mathcal
H_{\phi}).
\end{eqnarray}
Since $u_{\alpha\beta jj}=u_{ jj\alpha\beta}$, (\ref{2.31}) and
(\ref{2.33}) yield
\begin{eqnarray}\label{2.34}
&&F^{\alpha\beta}\phi_{\alpha\beta}- \phi_t\nonumber\\
&\leq& \sum_{j\in
B}u_n^{-3}\left[\sigma_l(G)+\frac{{\sigma}^2_1(B|j)-{\sigma}_2(B|j)}{{\sigma}^2_1(B)}\right]\left\{
\Big[\sum_{\alpha,\beta,r,s=1}^n
F^{\alpha\beta, rs}u_{\alpha\beta j}u_{rsj}\right.\nonumber\\
&&\qquad \qquad \qquad \qquad
+2\sum_{\alpha,\beta=1}^nF^{\alpha\beta,
u_n}u_{j\alpha\beta}u_{jn}+F^{u_n,u_n}u_{jn}^2\Big]u_n^2\nonumber\\
&&\qquad \qquad\qquad \qquad
+\left.4u_{jn}u_n\sum_{\alpha,\beta=1}^nF^{\alpha\beta}u_{j\alpha\beta}
-6u_{jn}^2\sum_{\alpha,\beta=1}^nF^{\alpha\beta}u_{\alpha\beta}\right\}\nonumber\\
&&-\frac{1}{{\sigma}^3_1(B)}\sum_{\alpha,\beta=1}^n\sum_{i\in
B}F^{\alpha\beta}[{\sigma}_1(B)a_{ii,\alpha}-a_{ii}\sum_{j\in
B}a_{jj,\alpha}][{\sigma}_1(B)a_{ii,\beta}-a_{ii}\sum_{j\in
B}a_{jj,\beta}]\nonumber\\
&&-\frac{1}{{\sigma}_1(B)}\sum_{\alpha,\beta=1}^n\sum_{i\neq j,
i,j\in B}F^{\alpha\beta}a_{ij,\alpha}a_{ij,\beta}+O(\mathcal H_{\phi}).
\end{eqnarray}
For each $ j \in B$, set
\begin{eqnarray}\label{2.35}
S_j= && \Big[\sum_{\alpha,\beta,r,s=1}^n
F^{\alpha\beta, rs}u_{j\alpha\beta}u_{rsj}+2\sum_{\alpha,\beta=1}^nF^{\alpha\beta,
u_n}u_{j\alpha\beta}u_{jn}+F^{u_n,u_n}u_{jn}^2\Big]u_n^2\nonumber\\
&+&4\sum_{\alpha,\beta=1}^nF^{\alpha\beta}u_{j\alpha\beta}u_{jn}u_n
-6\sum_{\alpha,\beta=1}^nF^{\alpha\beta}u_{\alpha\beta}u_{jn}^2
\end{eqnarray}
For each $j\in B$, set
\begin{eqnarray}\label{4.16}
&&X_{\alpha\beta}=u_{\alpha\beta j}u_n, \forall (\alpha, \beta), \\
&&Y=u_{jn}u_n.
\end{eqnarray}
In the coordinate system (\ref{2.6}),
\[(\n^2u(x), \n u(x), u(x), t)=(\n^2u,(0,...,0,|\n u|), u, t).\]
Equalize it to $(s^2A, s\theta, u, t)$, the components of
$\widetilde{V}$ defined in Lemma~\ref{lem2.2} are
\begin{eqnarray*}
&&\widetilde{X}_{\alpha\beta}=\frac{u_{\alpha\beta
j}}{u_n}-\frac{2u_{\alpha\beta}u_{jn}}{u^2_n}, \quad \forall
(\alpha,\beta),  \\
&&\widetilde{Y}=u_{jn}u_n.
\end{eqnarray*}
For $j \in B$, Lemma~\ref{lem2.2} implies
\begin{equation}\label{2.38}
S_j\le 0.
\end{equation}
Condition (\ref{1.2}) implies
\begin{equation}\label{2.39}
 (F^{\alpha\beta}) \ge \delta_0 I, \; \quad \mbox{for some $\delta_0>0$, and
$\forall x \in \O$.}
\end{equation}
Set
\[V_{i\alpha}={\sigma}_1(B)a_{ii,\alpha}-a_{ii}\sum_{j\in
B}a_{jj,\alpha}.\]
Combining (\ref{2.34}), (\ref{2.38}) and (\ref{2.39}),
\begin{eqnarray}\label{2.40}
F^{\alpha\beta}\phi_{\alpha\beta}\le  C(\phi+\sum_{i,j\in B}|\nabla
a_{ij}|)-\delta_0[\frac{\sum_{i\neq j\in
B, \alpha=1}^na^2_{ij\alpha}}{\sigma_1(B)}+\frac{\sum_{i\in
B,\alpha=1}^nV_{i\alpha}^2}{{\sigma}^3_1(B)}].
\end{eqnarray}
By Lemma~3.3 in \cite{BG09}, for each $M\ge 1$, for any $M\ge
|\gamma_i|\ge \frac{1}{M}$, there is a constant $C$ depending only
on $n$ and $M$ such that, $\forall \alpha$,
\begin{equation}\label{2.41}
\sum_{i,j\in B}|a_{ij\alpha}|\le C(1+\frac{1}{\delta_0^2})
(\sigma_1(B)+|\sum_{i\in B} \gamma_ia_{ii\alpha}|)
+\frac{\delta_0}{2}[\frac{\sum_{i\neq j \in B}|a_{ij\alpha}|^{2}}
{\sigma_{1}(B)}+
\frac{\sum_{i\in B}V_{i\alpha}^{2}}{\sigma_{1}^{3}(B)}].\end{equation}
Taking $\gamma_i=\sigma_l(G)+\frac{{\sigma}^2_1(B|i)-{\sigma}_2(B|i)}{{\sigma}^2_1(B)}$ for each $i\in B$,
the Newton-MacLaurine inequality implies
\[\sigma_l(G)+1 \ge
\sigma_l(G)+\frac{{\sigma}^2_1(B|j)-{\sigma}_2(B|j)}{{\sigma}^2_1(B)}
\ge \sigma_l(G), \quad \forall j\in B.\]
Therefore we conclude from Lemma~\ref{lem2.1} and (\ref{2.41}) that
$\sum_{i,j \in B} |\nabla a_{ij}|$  can be controlled by the rest
terms on the right hand side in (\ref{2.40}) and $\phi + |\nabla
\phi|$. The proof is complete. \qed

\section{Proof of Theorem \ref{1.2}}

In this section, through modifying the proof of Theorem \ref{1.1},
we will give a proof of Theorem \ref{1.2}. Also it is a parabolic
equation case corresponding to \cite{GX}.

Suppose that $u(x,t)$ is a spatial quasiconcave solution of \eqref{1.4},
and assume that level surface $\Sigma^{u(x_0,t)} = \{x \in \Omega|
u(x,t)=u(x_0,t)\}$ is connected for each $(x_0,t) \in \mathcal
{O}\times [0, T]$.

Set
\begin{equation} \label{3.1}
\widetilde a = a - \eta _0 g I, \quad \eta _0 \geqslant 0, \quad
g(x,t) = e^{Au(x,t)},
\end{equation}
where $A > 0$ is a constant to be determined. We want to show $\widetilde{a}$ is
of constant rank. Theorem \ref{1.1} corresponds to the case
$\eta_0=0$. If $ \min \{ \kappa ^0 ,\kappa ^1 \}=0$, there is
nothing to prove instead of utilizing Theorem \ref{1.1}. We will
assume $ \min \{ \kappa ^0 ,\kappa ^1 \}>0$ in the rest of the
paper. Denote $\kappa_s(x,t)$ and $\widetilde \kappa_s(x,t)$ be the
minimum eigenvalue of matrix $a(x)$ and $\widetilde a(x)$
respectively. Since the spatial level sets are strictly convex, and
$\overline{\Omega}$ is compact, $a$ is strictly positive definite.
That is, $\kappa_s(x,t)$ has a positive lower bound.

For a positive constant $A$ to be determined, increasing $\eta_0$
from 0, such that $\widetilde a$ is degenerate at some points, i.e.
$\widetilde a$ is semi-positive with the rank is not full. \eqref{1.8}
follows easily if this happens only on the boundary. We want to show that, if the degeneracy happens
at an interior point of $\Omega$,  then $\widetilde a$ is degenerate through out $\Omega$
 with the same rank. This
implies that the "=" holds in \eqref{1.8} and Theorem \ref{th1.2} is proved.

Therefore, the main task is to prove constant rank theorem for
$\widetilde a$. Suppose $\widetilde a(x,t_0)$ attains minimal rank
$l=l(t_0)$ at some point $z_0 \in \Omega$. We may assume $l\leqslant
n-2$, otherwise there is nothing to prove. And we assume $u \in
C^{3,1}$ and $u_n>0$ in the rest of this paper. So there is
a neighborhood $\mathcal {O}\times (t_0-\delta, t_0+\delta]$ of
$(z_0, t_0)$, such that there are $l$ "good" eigenvalues of
$(\widetilde a_{ij})$ which are bounded below by a positive
constant, and the other $n-1-l$ "bad" eigenvalues of $(\widetilde
a_{ij})$ are very small. Denote $G$ be the index set of these "good"
eigenvalues and $B$ be the index set of "bad" eigenvalues. And for
any fixed point $(x,t) \in \mathcal {O}\times (t_0-\delta,
t_0+\delta]$, we may express $(\widetilde a_{ij})$ in a form
of \eqref{3.1} and (\ref{2.5}), by choosing $e_1,\cdots,
e_{n-1},e_n$ such that
\begin{equation}\label{3.2}
|\n u(x,t)|=u_n(x,t)>0\ \mbox{and}\
(u_{ij}),i,j=1,..,n-1,\ \mbox{is diagonal at}\ (x,t).
\end{equation}

Without loss of generality, we assume $ u_{11} \geq u_{22}\geq \cdots
\geq u_{n-1,n-1} $. So, at $(x,t) \in \mathcal {O}\times (t_0-\delta,
t_0+\delta)$, from \eqref{2.5}, we have the matrix
$(a_{ij}),i,j=1,..,n-1$, is also diagonal. And without loss of
generality we may assume $a_{11} \geq a_{22} \geq ... \geq a_{n-1,
n-1}$, then $\widetilde a_{11} \geq \widetilde a_{22} \geq ... \geq
\widetilde a_{n-1, n-1}$. There is a positive constant $C>0$
depending only on $\|u\|_{C^{4}}$ and $\mathcal {O}\times
(t_0-\delta, t_0+\delta]$, such that $\widetilde a_{11} \ge
\widetilde a_{22} \ge...\ge \widetilde a_{ll} > C$ for all $(x,t)
\in \mathcal {O}\times (t_0-\delta, t_0+\delta)$. For convenience we
denote  $ G = \{ 1, \cdots ,l \} $ and $ B = \{
l+1, \cdots, n-1 \} $ be the "good" and "bad" sets of indices
respectively. If there is no confusion, we also denote
\begin{eqnarray}\label{3.3}
\quad \quad \mbox{\it $G=\{\widetilde a_{11},...,\widetilde
a_{ll}\}$ and $B=\{\widetilde a_{l+1,l+1},...,\widetilde a_{n-1,
n-1}\}$.}
\end{eqnarray}
Note that for any $\delta>0$, we may choose $\mathcal {O}\times
(t_0-\delta, t_0+\delta]$ small enough such that $\widetilde a_{jj}
<\delta$ for all $j \in B$ and $(x,t) \in \mathcal {O}\times
(t_0-\delta, t_0+\delta]$.

For each $(x,t)$, let $a=(a_{ij})$ be the symmetric Weingarten
tensor of $\Sigma^{u(x,t)}$. Set
\begin{eqnarray}\label{3.4}
p(\widetilde a)=\sigma_{l+1}(\widetilde a_{ij}), \quad q(\widetilde
a) &=& \left\{
\begin{array}{llr}
\frac{\sigma_{l+2}(\widetilde a_{ij})}{\sigma_{l+1}(\widetilde a_{ij})},  &\mbox{if} \; \sigma_{l+1}(\widetilde a_{ij})>0,&\\
0,& \mbox{otherwise}. &
\end{array}
\right.
\end{eqnarray}
Theorem \ref{th1.2} is equivalent to say $p(\widetilde a)\equiv 0$
(defined in (\ref{3.4}) ) in $\mathcal {O}\times (t_0-\delta, t_0]$.
As in the description of the proof of Theorem \ref{th1.1}, we should consider the function
\begin{equation}\label{3.5}
\phi(\widetilde a)= p(\widetilde a)+q(\widetilde a),
\end{equation}
where $p$ and $q$ as in \eqref{3.4}. We will use notion $h=O(f)$ if
$|h(x,t)| \le Cf(x,t)$ for $(x,t) \in \O \times(t_0-\delta, t_0+\delta]$ with positive constant $C$ under
control.

To get around $p=0$, for $\varepsilon>0$ sufficiently small, we
instead consider
\begin{equation}\label{3.6}
\phi_\varepsilon (\widetilde a)= \phi(\widetilde a_\varepsilon),
\end{equation}
where $a_\varepsilon=\widetilde a+\varepsilon I.$ We will also
denote $G_\varepsilon=\{\widetilde a_{ii}+\varepsilon, i\in G\},$
$B_\varepsilon=\{\widetilde a_{ii}+\varepsilon, i\in B\}.$

To simplify the notations, we will drop subindex $\varepsilon$ with
the understanding that all the estimates will be independent of
$\varepsilon.$ In this setting, if we pick $\mathcal {O}\times
(t_0-\delta, t_0+\delta]$ small enough, there is $C>0$ independent
of $\varepsilon$ such that
\begin{equation}\label{2.11}
\phi(\widetilde a(x,t))\geq C\varepsilon, \quad \sigma_1(B)\geq
C\varepsilon, ~\quad ~\rm{ for ~all~ }(x,t) \in \mathcal {O}\times
(t_0-\delta, t_0+\delta].
 \end{equation}
We also denote
\begin{equation} \mathcal H_{\phi} = \sum_{i,j\in B}|\nabla \widetilde a_{ij}|+\phi.
\end{equation}

\begin{lemma}\label{lem3.3}
For any fixed $(x,t) \in \mathcal {O}\times (t_0-\delta,
t_0+\delta]$, with the coordinate chosen as in (\ref{3.2}) and
(\ref{3.3}),

\begin{eqnarray}\label{3.9}
&&\sum_{\alpha,\beta=1}^{n}F^{\alpha\beta}\phi_{\alpha\beta}-\phi_t
\notag  \\
&=&u_n^{-3}\sum_{j\in B}\left[ \sigma_l(G) +
\frac{{\sigma}^2_1(B|j)-{\sigma}_2(B|j)} {{\sigma}^2_1(B)}\right]
[-u_n^2(\sum_{\alpha,\beta=1}^n F^{\alpha\beta}u_{\alpha\beta j
j}-u_{jjt})+2u_nu_{nj}(\sum_{\alpha,\beta=1}^nF^{\alpha\beta}u_{\alpha\beta
j}-u_{jt})
\notag \\
&& \qquad\qquad\qquad\qquad\qquad\qquad\qquad\qquad\quad
+4u_nu_{nj}\sum_{\alpha,\beta=1}^nF^{\alpha\beta}u_{\alpha\beta
j}-6u_{nj}^2\sum_{\alpha,\beta=1}^{n}F^{\alpha\beta}u_{\alpha\beta
}]  \notag  \\
&&+2u_n^{-3}\sum_{j\in B,i\in G}\left[ \sigma_l(G) +
\frac{{\sigma}^2_1(B|j)-{\sigma}_2(B|j)}
{{\sigma}^2_1(B)}\right]\sum_{\alpha,\beta=1}^{n}
F^{\alpha\beta}\frac{1}{u_{ii}}[u_nu_{ij\alpha}-2u_{i\a}u_{jn}][u_nu_{ij\b}-2u_{i\b}u_{jn}] \notag \\
&&+\eta _0 g \left[-A^2F^{nn}u_n^2+A O(1) + O (1)\right] \notag   \\
&&-\frac{1}{{\sigma}^3_1(B)}\sum_{\alpha,\beta=1}^n\sum_{i\in
B}F^{\alpha\beta}[{\sigma}_1(B)\widetilde a_{ii,\alpha}-\widetilde
a_{ii}\sum_{j\in B}\widetilde a_{jj,\alpha}][{\sigma}_1(B)\widetilde
a_{ii,\beta}-\widetilde a_{ii}\sum_{j\in B}\widetilde a_{jj,\beta}]\nonumber\\
&& -\frac{1}{{\sigma}_1(B)}\sum_{\alpha,\beta=1}^n\sum_{i\neq j\in
B}F^{\alpha\beta}\widetilde a_{ij,\alpha} \widetilde
a_{ij,\beta}+O(\mathcal H_{\phi}). \notag
\end{eqnarray}
\end{lemma}

{\bf Proof:} For any fixed $(x,t) \in \mathcal {O}\times (t_0-\delta,
t_0+\delta]$, we choose the coordinate  as in (\ref{3.2}) such that
$|\n u(x)|= u_n(x) >0$ and the matrix $(\widetilde a_{ij}(x))$ is
diagonal for $1 \le i,j\le n-1$ and nonnegative. From the definition
of $p$,
\begin{eqnarray}\label{3.9}
a_{jj}=-\frac{h_{jj}}{u^3_n}=-\frac{u_{jj}}{u_n}=O(\mathcal
H_{\phi}), \forall j\in B,
\end{eqnarray}
and
\begin{equation}\label{3.10}
\phi_t=\sum_{j \in B} \left[ \sigma_l(G) +
\frac{{\sigma}^2_1(B|j)-{\sigma}_2(B|j)} {{\sigma}^2_1(B)}\right]
\widetilde a_{jjt} +O(\mathcal H_{\phi}).
\end{equation}
Using relationship (\ref{3.9}), we have
\begin{eqnarray}\label{3.11}
\phi_{\alpha \beta}&=&\sum_{j\in B}\left[ \sigma_l(G) +
\frac{{\sigma}^2_1(B|j)-{\sigma}_2(B|j)}
{{\sigma}^2_1(B)}\right]\Big[\widetilde a_{jj,\alpha
\beta}-2\sum_{i\in G}\frac{\widetilde a_{ij,\alpha}\widetilde
a_{ij,\beta}}{\widetilde a_{ii}}\Big] \nonumber\\
&&-\frac{1}{{\sigma}^3_1(B)}\sum_{i\in
B}\Big[{\sigma}_1(B)\widetilde a_{ii,\alpha}-\widetilde
a_{ii}\sum_{j\in B}\widetilde a_{jj,\alpha} \Big] \Big[{\sigma}_1(B)
\widetilde a_{ii,\beta} - \widetilde a_{ii}\sum_{j\in
B} \widetilde a_{jj,\beta} \Big]  \nonumber\\
&&-\frac{1}{{\sigma}_1(B)}\sum_{i\neq j\in B}\widetilde
a_{ij,\alpha}\widetilde a_{ij,\beta}+O(\mathcal H_{\phi}).
\end{eqnarray}

So,
\begin{eqnarray}\label{3.12}
&&\sum_{\alpha,\beta=1}^{n}F^{\alpha\beta}[\widetilde a_{jj,\alpha
\beta}-2\sum_{i\in G}\frac{\widetilde a_{ij,\alpha}\widetilde
a_{ij,\beta}}{\widetilde a_{ii}}] -\widetilde a_{jj,t} \notag  \\
&=&\sum_{\alpha,\beta=1}^{n}F^{\alpha\beta}a_{jj,\alpha
\beta}-a_{jj,t}+\sum_{\alpha,\beta=1}^{n}F^{\alpha\beta}(-\eta_0g_{\a\b})+\eta_0g_t
\notag \\
&&-2\sum_{\alpha,\beta=1}^{n}F^{\alpha\beta}\sum_{i\in
G}\frac{\widetilde a_{ij,\alpha}\widetilde a_{ij,\beta}}{\widetilde
a_{ii}}.
\end{eqnarray}
From the definition of $a_{ij}$, and $u_k = 0$ for $k = 1, \cdots
, n - 1$, we can get
\begin{eqnarray}\label{3.13}
u_n u_{ij\alpha }  =  - u_n ^2 a_{ij,\alpha }  + u_{nj} u_{i\alpha }
+ u_{ni} u_{j\alpha }  + u_{n\alpha } u_{ij}
\end{eqnarray}
and
\begin{eqnarray}\label{3.14}
u_n ^3 a_{jj,\alpha \beta } & = & - u_n ^2 u_{jj\alpha \beta }  +
2u_n u_{nj} u_{\alpha \beta j}  - 2u_n (u_{n\beta } u_{jj\alpha }  +
u_{n\alpha } u_{jj\beta } )   \notag \\
&& + 2u_n (u_{j\alpha } u_{nj\beta }  + u_{j\beta } u_{nj\alpha } )
+ 2u_{nj} (u_{n\alpha } u_{j\beta }  + u_{n\beta } u_{j\alpha } ) -
2u_{nn} u_{j\alpha } u_{j\beta }   \notag \\
&& - 2(u_{n\alpha } u_{n\beta }  + u_n u_{\alpha \beta n} )u_{jj}  -
2\eta _0 gu_{j\alpha } u_{j\beta } u_n  - 3\eta _0 u_n ^2
(u_{n\alpha } g_\beta   + u_{n\beta } g_\alpha  )  \notag \\
&& - \eta _0 g(3u_n ^2 u_{n\alpha \beta }  + 6u_{n\alpha } u_{n\beta
} u_n  + \sum\limits_{i = 1}^{n - 1} {u_{i\alpha } u_{i\beta } u_n }
) + O(\mathcal H_{\phi}).
\end{eqnarray}
Direct calculation and \eqref{3.13}, we can get
\begin{eqnarray}\label{3.15}
&&- a_{jj,t}  + \sum\limits_{\alpha \beta  = 1}^n {F^{\alpha \beta }
( - \eta _0 g_{\alpha \beta } )}  + \eta _0 g_t  \notag  \\
& =& \frac{1} {{u_n ^3 }}[u_n ^2 u_{jjt}  - 2u_n u_{nj} u_{jt} ]
\notag  \\
&&+ \eta _0 g[ - A^2 F^{nn} u_n ^2  - A(\sum\limits_{\alpha \beta  =
1}^n {F^{\alpha \beta } u_{\alpha \beta } }  - u_t ) + \frac{{u_{nt}
}} {{u_n }}].
\end{eqnarray}
From \eqref{3.14},
\begin{eqnarray}
\sum\limits_{\alpha \beta  = 1}^n {F^{\alpha \beta } a_{jj,\alpha
\beta } }  &=& \sum\limits_{\alpha \beta  = 1}^n {\frac{{F^{\alpha
\beta } }} {{u_n ^3 }}[ - u_n ^2 u_{jj\alpha \beta }  + 2u_n
u_{nj}u_{\alpha \beta j} }   \notag  \\
&& - 4u_{nj} u_{n\alpha } u_{j\beta }  + 4u_{nj} u_{n\alpha }
u_{j\beta }  - 2u_{nn} u_{j\alpha } u_{j\beta }    \notag \\
&& - 2\eta _0 u_n ^2 u_{n\alpha } g_\beta   - \eta _0 g(u_n ^2
u_{n\alpha \beta }  + 2u_{j\alpha } u_{j\beta } u_n  +
\sum\limits_{i = 1}^{n - 1} {u_{i\alpha } u_{i\beta } u_n } )]
+O(\mathcal H_{\phi}),   \notag
\end{eqnarray}
so, as in \cite{GX}, we can get
\begin{eqnarray} \label{3.16}
&&u_n ^3 \sum\limits_{\alpha \beta  = 1}^n {F^{\alpha \beta }
a_{jj,\alpha \beta } }   \notag  \\
& =&  - u_n ^2\sum\limits_{\alpha \beta = 1}^n {F^{\alpha \beta }
u_{jj\alpha \beta } }  + 2u_n u_{nj} \sum\limits_{\alpha \beta  =
1}^n {F^{\alpha \beta } u_{\alpha \beta
j} }   \notag   \\
&& + 4u_n u_{nj} \sum\limits_{\alpha \beta  = 1}^n {F^{\alpha \beta
} u_{\alpha \beta j} }  - 6u_{nj} ^2 \sum\limits_{\alpha \beta  =
1}^n {F^{\alpha \beta } u_{\alpha \beta } }   \notag  \\
&& + 4u_n ^2 \sum\limits_{\alpha  = 1}^n {\sum\limits_{i \in G}
{F^{\alpha i} a_{ij,\alpha } } }  + 2u_{nj} ^2 \sum\limits_{i \in G}
{F^{ii} u_{ii} }   \notag  \\
&& + 2u_{nj} ^2 \sum\limits_{i \in B} {F^{ii} u_{ii} }  - 12u_{jn}
u_{jj} \sum\limits_{\alpha  = 1}^n {F^{j\alpha } u_{n\alpha } }  +
4u_n u_{jj} \sum\limits_{\alpha  = 1}^n {F^{j\alpha } u_{jn\alpha }
}  - 2u_{nn} F^{jj} u_{jj} ^2  \notag   \\
&& - \eta _0 g(u_n ^2 u_{n\alpha \beta }  + 2u_{j\alpha } u_{j\beta
} u_n  + \sum\limits_{i = 1}^{n - 1} {u_{i\alpha } u_{i\beta } u_n }
)   \notag \\
&& - 2\eta _0 \sum\limits_{\alpha \beta  = 1}^n {F^{\alpha \beta }
u_{n\alpha } g_\beta  } u_n  + 4\eta _0 \sum\limits_{\alpha  = 1}^n
{F^{j\alpha } g_\alpha  u_{jn} u_n ^2 }  + O(\mathcal H_{\phi}),
\notag   \\
& =&  - u_n ^2\sum\limits_{\alpha \beta = 1}^n {F^{\alpha \beta }
u_{jj\alpha \beta } }  + 2u_n u_{nj} \sum\limits_{\alpha \beta  =
1}^n {F^{\alpha \beta } u_{\alpha \beta j} }   \notag  \\
&& + 4u_n u_{nj} \sum\limits_{\alpha \beta  = 1}^n {F^{\alpha \beta
} u_{\alpha \beta j} }  - 6u_{nj} ^2 \sum\limits_{\alpha \beta  =
1}^n {F^{\alpha \beta } u_{\alpha \beta } }     \\
&& + 4u_n ^2 \sum\limits_{\alpha  = 1}^n {\sum\limits_{i \in G}
{F^{\alpha i} a_{ij,\alpha } } }  + 2u_{nj} ^2 \sum\limits_{i \in G}
{F^{ii} u_{ii} }   \notag  \\
&& +\eta _0 g \left[A O(1) + O (1)\right]  + O(\mathcal H_{\phi}).
\notag
\end{eqnarray}
Also, with the similar computations \eqref{2.23} in the Lemma \ref{lem2.1},
\begin{eqnarray} \label{3.17}
&&\sum\limits_{\alpha \beta  = 1}^n {F^{\alpha \beta }
\frac{{\widetilde a_{ij,\alpha } \widetilde a_{ij,\beta } }}
{{\widetilde a_{ii} }}}  - \frac{1} {{u_n ^3 }}[2u_n ^2 u_{nj}
\sum\limits_{\alpha  = 1}^n {F^{\alpha i} a_{ij,\alpha } }  + u_{nj}
^2 F^{ii} u_{ii} ]   \notag  \\
& =& \sum\limits_{\alpha \beta  = 1}^n {F^{\alpha \beta }
\frac{{a_{ij,\alpha } a_{ij,\beta } }} {{a_{ii} }}}  + \eta _0
g\sum\limits_{\alpha \beta  = 1}^n {F^{\alpha \beta }
\frac{{a_{ij,\alpha } a_{ij,\beta } }} {{a_{ii} \widetilde a_{ii}
}}}   \notag  \\
&&- \frac{1} {{u_n ^3 }}[2u_n ^2 u_{nj} \sum\limits_{\alpha  = 1}^n
{F^{\alpha i} a_{ij,\alpha } }  + u_{nj} ^2 F^{ii} u_{ii} ]  \notag \\
& =& \eta _0 g\sum\limits_{\alpha \beta  = 1}^n {F^{\alpha \beta }
\frac{{a_{ij,\alpha } a_{ij,\beta } }} {{a_{ii} \widetilde a_{ii}
}}}  \notag  \\
&&- \frac{1} {{u_n ^3 }}\sum\limits_{\alpha \beta  = 1}^n {F^{\alpha
\beta } \frac{1} {{u_{ii} }}} [ - u_n u_{ij\alpha }  + u_{nj}
u_{i\alpha } + u_{ni} u_{j\alpha } ][ - u_n u_{ij\beta }  + u_{nj}
u_{i\beta }  + u_{ni} u_{j\beta } ]   \notag  \\
&& - \frac{1} {{u_n ^3 }}[2u_{nj} \sum\limits_{\alpha  = 1}^n
{F^{\alpha i} ( - u_n u_{ij\alpha }  + u_{nj} u_{i\alpha }  + u_{ni}
u_{j\alpha } )} + u_{nj} ^2 F^{ii} u_{ii} ]   \notag  \\
& =& \eta _0 g\sum\limits_{\alpha \beta  = 1}^n {F^{\alpha \beta }
\frac{{a_{ij,\alpha } a_{ij,\beta } }} {{a_{ii} \widetilde a_{ii}
}}}  \notag  \\
&& - \frac{1} {{u_n ^3 }}\sum\limits_{\alpha \beta  = 1}^n
{F^{\alpha \beta } \frac{1} {{u_{ii} }}} [ - u_n u_{ij\alpha }  +
2u_{nj} u_{i\alpha } ][ - u_n u_{ij\beta }  + 2u_{nj} u_{i\beta } ]
\notag  \\
&& - \frac{1} {{u_n ^3 }}u_{jj} [\sum\limits_{\alpha  = 1}^{n - 1}
{F^{\alpha j} \frac{2} {{u_{ii} }}} u_{ni} ( - u_n u_{ij\alpha }  +
u_{nj} u_{i\alpha } ) + F^{ii} \frac{1} {{u_{ii} }}u_{jj} u_{ni} ^2
\notag  \\
&& \qquad \qquad + 2F^{jn} \frac{1} {{u_{ii} }}u_{ni} ( - u_n
u_{ijn}  + 2u_{nj} u_{in} ) + F^{ij} u_{ni} u_{nj} ]   \notag  \\
& =&  - \frac{1} {{u_n ^3 }}\sum\limits_{\alpha \beta  = 1}^n
{F^{\alpha \beta } \frac{1} {{u_{ii} }}} [ - u_n u_{ij\alpha }  +
2u_{nj} u_{i\alpha } ][ - u_n u_{ij\beta }  + 2u_{nj} u_{i\beta } ]
\notag   \\
&& + \eta _0 g\left[\sum\limits_{\alpha \beta  = 1}^n {F^{\alpha
\beta } \frac{{a_{ij,\alpha } a_{ij,\beta } }} {{a_{ii} \widetilde
a_{ii}}}} +O(1). \right]
\end{eqnarray}

From the above calculations, the proof is complete.  \qed

Theorem~\ref{th1.2} is a direct consequence of the following
proposition and the strong maximum principle.

\begin{proposition}\label{thm4.1}
Suppose that the function $F, u$ satisfy assumptions in Theorem
\ref{th1.2}. If the second fundamental form $b_{ij}$ of
$\Sigma^{u(x,t_0)}$ attains minimum rank $l=l(t_0)$ at certain point
$x_0\in \Omega$, then there exist a neighborhood $\mathcal {O}
\times(t_0-\delta_0,t_0+\delta_0]$ of $(x_0,t_0)$ and a positive
constant $C$ independent of $\phi$ (defined in (\ref{3.5})), such
that
\begin{equation}\label{3.18}
\sum_{\alpha,\beta=1}^{n}F^{\alpha\beta}\phi_{\alpha\beta}(x,t)-\phi_t\leq
C(\phi+|\nabla \phi|)+\eta _0 g \left[-A^2F^{nn}u_n^2+A O(1) + O
(1)\right]
\end{equation}
holds for any $(x,t)\in \mathcal {O} \times (t_0-\delta_0,
t_0+\delta_0]$.
\end{proposition}

{\bf Proof:} Since
\begin{equation}\label{3.19}
u_{t}=F(\n ^2 u, \n u ,u, t),
\end{equation}
for each $j\in B$, differentiating the above equation in $e_j$
direction at $x$,
\begin{equation}\label{3.20}
u_{jt}=\sum_{\alpha,\beta=1}^nF^{\alpha\beta}u_{\alpha\beta
j}+F^{u_n}u_{jn}+O(\mathcal H_{\phi})
\end{equation}
and
\begin{eqnarray}\label{3.21}
u_{jjt}&=&\sum_{\alpha,\beta=1}^nF^{\alpha\beta}u_{\alpha\beta
jj}+\sum_{\alpha,\beta,r,s=1}^nF^{\alpha\beta, rs}u_{\alpha\beta
j}u_{rsj} +2\sum_{\alpha,\beta,l=1}^nF^{\alpha\beta,
u_l}u_{\alpha\beta j}u_{lj}\nonumber\\
&&+2\sum_{\alpha,\beta=1}^nF^{\alpha\beta, u}u_{j\alpha\beta}u_j
+\sum_{l,s=1}^nF^{u_l, u_s}u_{lj}u_{sj}-2\sum_{l=1}^nF^{u_l,
u}u_{lj}u_j\nonumber\\
&&+F^{u,u}u_j^2+\sum_{l=1}^nF^{u_l}u_{ljj}+F^{u}u_{jj}.
\end{eqnarray}
It follows from (\ref{3.9})and \eqref{3.13} that, at $(x,t)$
\begin{eqnarray}
 \sum\limits_{\alpha \beta  = 1}^n {F^{\alpha \beta } u_{\alpha
\beta j} }  - u_{jt}  =  - F^{p_n } u_{nj}  + \eta _0 gF^{p_j } u_n
+ O(\mathcal H_{\phi})
\end{eqnarray}
and
\begin{eqnarray}
&&\sum\limits_{\alpha \beta  = 1}^n {F^{\alpha \beta } u_{\alpha
\beta jj} }  - u_{jjt}  \notag  \\
&=&  - \sum\limits_{\alpha \beta \gamma \eta = 1}^n {F^{\alpha \beta
,\gamma \eta } u_{\alpha \beta j} } u_{\gamma \eta j}  -
2\sum\limits_{\alpha \beta  = 1}^n {F^{\alpha \beta ,p_n } u_{\alpha
\beta j} u_{nj} } - F^{p_n ,p_n } u_{nj} u_{nj}  - 2\frac{{F^{p_n }
}} {{u_n }}u_{nj} ^2   \notag  \\
&& + \eta _0 g[ - AF^{p_n } u_n ^2 ]   \notag  \\
&& + \eta _0 g[2\sum\limits_{\alpha \beta  = 1}^n {F^{\alpha \beta
,p_j } u_{\alpha \beta j} u_n }  + F^{p_j ,p_j } u_{jj} u_n  +
2F^{p_n ,p_j } u_{nj} u_n  + F^{p_n } u_n  + 2F^{p_j } u_{jn}  +
F^{p_l } u_{nl} ]  \notag  \\
&&+ O(\mathcal H_{\phi}).   \notag
\end{eqnarray}
From lemma \ref{3.1},
\begin{eqnarray}\label{4.13}
&&F^{\alpha\beta}\phi_{\alpha\beta}- \phi_t\nonumber\\
&=& \sum_{j\in
B}u_n^{-3}\left[\sigma_l(G)+\frac{{\sigma}^2_1(B|j)-{\sigma}_2(B|j)}{{\sigma}^2_1(B)}\right]\left\{
\Big[\sum_{\alpha,\beta,r,s=1}^n
F^{\alpha\beta, rs}u_{\alpha\beta j}u_{rsj}\right.\nonumber\\
&&\qquad \qquad \qquad \qquad
+2\sum_{\alpha,\beta=1}^nF^{\alpha\beta,
u_n}u_{j\alpha\beta}u_{jn}+F^{u_n,u_n}u_{jn}^2\Big]u_n^2\nonumber\\
&&\qquad \qquad\qquad \qquad
+\left.4u_{jn}u_n\sum_{\alpha,\beta=1}^nF^{\alpha\beta}u_{j\alpha\beta}
-6u_{jn}^2\sum_{\alpha,\beta=1}^nF^{\alpha\beta}u_{\alpha\beta}\right\}\nonumber\\
&&+2u_n^{-3}\sum_{j\in B,i\in G}\left[ \sigma_l(G) +
\frac{{\sigma}^2_1(B|j)-{\sigma}_2(B|j)}
{{\sigma}^2_1(B)}\right]\sum_{\alpha,\beta=1}^{n}
F^{\alpha\beta}\frac{1}{u_{ii}}[u_nu_{ij\alpha}-2u_{i\a}u_{jn}][u_nu_{ij\b}-2u_{i\b}u_{jn}] \notag \\
&&+\eta _0 g \left[-A^2F^{nn}u_n^2+A O(1) + O (1)\right] \notag   \\
&&-\frac{1}{{\sigma}^3_1(B)}\sum_{\alpha,\beta=1}^n\sum_{i\in
B}F^{\alpha\beta}[{\sigma}_1(B)\widetilde a_{ii,\alpha}-\widetilde
a_{ii}\sum_{j\in B}\widetilde a_{jj,\alpha}][{\sigma}_1(B)\widetilde
a_{ii,\beta}-\widetilde a_{ii}\sum_{j\in B}\widetilde a_{jj,\beta}]\nonumber\\
&& -\frac{1}{{\sigma}_1(B)}\sum_{\alpha,\beta=1}^n\sum_{i\neq j\in
B}F^{\alpha\beta}\widetilde a_{ij,\alpha} \widetilde
a_{ij,\beta}+O(\mathcal H_{\phi}). \notag
\end{eqnarray}
So, following the argument in the proof of Proposition \ref{thm2.3},  we get,
\begin{equation}
\sum_{\alpha,\beta=1}^{n}F^{\alpha\beta}\phi_{\alpha\beta}(x,t)-\phi_t\leq
C(\phi+|\nabla \phi|)+\eta _0 g \left[-A^2F^{nn}u_n^2+A O(1) + O
(1)\right].
\end{equation}
The proof is completed. \qed

\bibliographystyle{amsplain}

\begin{thebibliography}{10}

\bibitem {AH} Ahlfors, L.V.: \textit{Conformal invariants: topics in geometric function theory}. McGraw-Hill Series in Higher Mathematics, McGraw-Hill Book Co., New York-D¨¹sseldorf-Johannesburg(1973)

\bibitem {BG09} Bian, B., Guan, P.:  A microscopic convexity principle for nonlinear partial differential equations.
Inventiones Math. {\bf 177}, 307-335(2009)

\bibitem{BGMX} Bian, B., Guan, P., Ma, X.N.,  Xu, L.: A constant rank theorem for quasiconcave solutions of fully nonlinear partial differential
equations. to appear in Indiana Univ. Math. J..

\bibitem {BLS} Bianchini, C., Longinetti, M., Salani, P.: Quasiconcave solutions to elliptic problems in convex rings, Indiana Univ. Math. J. {\bf 58}, 1565-1590(2009)

\bibitem {Bo} Borell, C.:  Brownian motion in a convex ring and quasi-concavity. Commun. Math. Phys. {\bf 86}, 143-147(1982)

\bibitem {CF85} Caffarelli, L., Friedman, A.: Convexity of solutions of some semilinear elliptic equations. Duke Math. J. {\bf 52}, 431--455(1985)

\bibitem{CS82} Caffarelli,  L., Spruck, J.: Convexity properties of solutions to some classical variational problems. Comm. Part. Diff. Eq. {\bf 7}, 1337-1379(1982)

\bibitem{CMY} Chang, S.-Y.A.,  Ma, X.N., Yang, P.: Principal curvature estimates for the convex level sets of semilinear elliptic equations. Discrete Contin. Dyn. Syst. {\bf 28}, 1151--1164(2010)

\bibitem{Ga57} Gabriel, R.: A result concerning convex level surfaces of 3-dimensional
harmonic functions.  J. London Math. Soc. {\bf 32}, 286-294(1957)

\bibitem{GM03} Guan, P., Ma, X.N.: The Christoffel-Minkowski problem I: convexity of solutions of a Hessian equations. Inventiones Math. {\bf 151}, 553-577(2003)

\bibitem{GX}Guan, P., Xu, L.:  Convexity estimates for level surfaces of quasiconcave solutions to fully nonlinear elliptic equations.
http://arxiv.org/abs/1004.1187v1

\bibitem{Ka85} Kawhol, B.:  \textit{Rearrangements and convexity of level sets in PDE}. Springer Lecture Notes in Math.{\bf 1150}(1985)

\bibitem{Ko90} Korevaar, N.: Convexity of level sets for solutions to elliptic ring problems. Comm. Part. Diff. Eq. {\bf 15}(4), 541-556(1990)

\bibitem{Le77} Lewis, J.: Capacitary functions in convex rings. Arch. Rat. Mech. Anal. {\bf 66}, 201-224(1977)

\bibitem{Lo83} Longinetti, M.: Convexity of the level lines of harmonic functions. (Italian) Boll. Un. Mat. Ital. A {\bf 6}, 71--75(1983)

\bibitem{Lo87} Longinetti, M.: On minimal surfaces bounded by two convex curves in parallel planes. J. Diff. Equations {\bf 67}, 344--358(1987)

\bibitem{MOZ09} Ma, X.N., Ou, Q.Z., Zhang, W.: Gaussian curvature estimates for the convex level sets of $p$-harmonic functions. Comm. Pure Appl. Math.{\bf 63}, 0935--0971(2010)
\bibitem{MZ} Ma, X.N., Zhang, W.:  The concavity of the Gaussian curvature of the convex level sets of $p$-harmonic functions with respect to the height. Preprint

\bibitem{OS} Ortel, M., Schneider, W.: Curvature of level curves of harmonic functions. Canad. Math. Bull. {\bf 26}, 399--405(1983)

\bibitem{Sh56} Shiffman, M.: On surfaces of stationary area bounded by two circles or convex curves in parallel planes. Annals of Math. {\bf 63}, 77--90(1956)

\bibitem{T83}  Talenti, G.: On functions whose lines of steepest descent bend proportionally to level lines. Ann. Scuola Norm. Sup. Pisa Cl. Sci. {\bf 10}(4), 587--605(1983)

\bibitem{WZ} Wang, P.H., Zhang, W.:  Gaussian curvature estimates for the convex level sets for some nonlinear partial differential equations. http://arxiv.org/abs/1003.2057v1

\bibitem{Xu08} Xu, L.: A  microscopic convexity theorem of level sets for solutions to elliptic equations. Cal. Var. PDE. DOI 10.1007/s00526-010-0333-3(2010)

\end{thebibliography}

\end{document}